\documentclass[11pt]{article}

\usepackage{latexsym,color,graphicx,amsmath,amssymb}

\usepackage{enumitem}

\def\reals{{\mathbb R}}
\def\C{{\cal C}}
\def\D{{\cal D}}

\newtheorem{theorem}{Theorem}[section]

\newtheorem{proposition}[theorem]{Proposition}
 
\begin{document}

\title{Distinct distances from three points\thanks{%
Work on this paper by Micha Sharir was 
supported by Grants 338/09 and 892/13 from the Israel Science Foundation,
by Grant 2012/229 from the US-Israel Binational Science Foundation,
by the Israeli Centers of Research Excellence (I-CORE)
program (Center No.~4/11), and
by the Hermann Minkowski-MINERVA Center for Geometry
at Tel Aviv University. 
Work by J\'ozsef Solymosi was supported by NSERC, 
ERC-AdG 321104, and OTKA NK 104183 grants. }}

\author{
Micha Sharir\thanks{%
School of Computer Science, Tel Aviv University,
Tel Aviv 69978, Israel.
{\sl michas@tau.ac.il} }
\and
J\'ozsef Solymosi\thanks{%
Department of Mathematics,
University of British Columbia,
Vancouver, BC, V6T 1Z4, Canada. 
{\sl solymosi@math.ubc.ca}} }

\maketitle

\begin{abstract}
Let $p_1,p_2,p_3$ be three noncollinear points in the plane, and let $P$ be 
a set of $n$ other points in the plane. We show that the number of distinct 
distances between $p_1,p_2,p_3$ and the points of $P$ is $\Omega(n^{6/11})$, 
improving the lower bound $\Omega(n^{0.502})$ of Elekes and Szab\'o~\cite{ESz} 
(and considerably simplifying the analysis).
\end{abstract}


%%%%%%%%%%%%%%%%%%%%%%%%%%%%%%%%%%%%%%%%%%%%%%%%%%%%%%%%%%%%%%%%%%%%%%%%%%%%%%%%

\section{Introduction}

The problem studied in this paper, as stated in the abstract, was raised by 
Erd\H os, Lov\'asz, and Vesztergombi \cite{ELV},
who conjectured that in the plane the number of distinct distances between three 
points, $p_1,p_2,p_3$, and $n$ other points is linear in $n$. This conjecture
was refuted by Elekes and Szab\'o \cite{ESz}, who gave a construction where
the number of distinct distances can be as small as $c\sqrt{n}$, for a suitable 
constant $c$, when $p_1,p_2$, and $p_3$ are collinear. Nevertheless, they also 
showed that if the three points are not collinear then there is a gap---the 
number of distinct distances is at least $n^{0.502}$. Using a different 
approach, which also appears to be considerably simpler, we improve this 
lower bound, for noncollinear $p_1,p_2,p_3$, to $\Omega(n^{6/11})$.

\paragraph{The general setup.}
Our derivation can be viewed as a special instance of a more general 
technique, which applies to the following general setup, as studied by 
Elekes and R\'onyai~\cite{ER00} and by Elekes and Szab\'o~\cite{ESz} 
(see also \cite{El}).
We have three sets $A$, $B$, $C$, each of $n$ real numbers, and we have
a trivariate real polynomial $F$ of degree $d$, which we assume to be
some constant. 
Let $Z(F)$ denote the subset of $A\times B\times C$ where $F$ vanishes.
Then, unless $F$ and $A$, $B$, $C$ have some very special structure,
$|Z(F)|$ should be subquadratic. 
(For a simple example where $|Z(F)|$ is quadratic in $n$, consider the 
case where $F(x,y,z)=p(x)+q(y)+r(z)$, for three suitable univariate 
polynomials $p$, $q$, and $r$, and where the respective images of 
$A$, $B$, and $C$ under $p$, $q$, and $r$ are, say, $\{1,2,\ldots,n\}$.)

Positive and significant results for this general problem have
been obtained by Elekes and R\'onyai~\cite{ER00} and by Elekes and 
Szab\'o~\cite{ESz}, who showed that if $|Z(F)| = \Omega(n^{1.95})$
and $n$ is large enough, then $F$ must indeed have a very restricted form.
For example, in the case where $F$ is of the form $z-f(x,y)$, $f$ must
be of the form $p(q(x)+r(y))$ or $p(q(x)\cdot r(y))$ for suitable polynomials
or rational functions $p$, $q$, $r$. Related representations, somewhat
more complicated to state, have also been obtained for the general case.

As will be apparent from the analysis in the following section, 
the problem that we study fits into this general scenario, 
for appropriate choices of $A$, $B$, $C$, and $F$. However, instead of
applying the general results reviewed above, we tackle the problem
in a more explicit and ad-hoc manner, which reduces the problem to an 
incidence problem between points and curves in a suitable parametric plane.

Our approach also applies to the general problem, and, in this context,
it can be (briefly) described as follows.\footnote{%
  The general technique, as described next, is incomplete, in the sense
  that we do not yet have a way to handle, in full generality, one crucial 
  step in the analysis (concerning the multiplicity of certain curves
  constructed by the analysis; see below). We provide this general approach
  to put our problem in the appropriate more general perspective, and to
  raise the open problem of closing this gap in the analysis.}
Let $A$, $B$, $C$, and $F$ be as above, and put $M=|Z(F)|$.
For each $a\in A$, $b\in B$, consider the planar curve
$\gamma_{a,b}$, which is the locus of all $(x,y) \in A\times B$
for which there exists $z\in C$ such that $F(x,b,z)=F(a,y,z)=0$.

Let $\Pi$ denote the set $A\times B$ in the $xy$-plane, let
$\Gamma$ denote the (multi-)set of the curves $\gamma_{a,b}$, 
and let $I=I(\Pi,\Gamma)$ denote the number of incidences 
between the curves of $\Gamma$ and the points of $\Pi$. 

For each $c\in C$, put
$$
\Pi_c = \{(x,y)\in \Pi \mid F(x,y,c)=0 \} ,
$$
and put $M_c = |\Pi_c|$.  We clearly have $\sum_{c\in C} M_c = M$.

Fix $c\in C$, and note that for any pair of pairs $(a_1,b_1)$, 
$(a_2,b_2)\in \Pi_c$, we have $(a_1,b_2)\in\gamma_{a_2,b_1}$
and $(a_2,b_1)\in\gamma_{a_1,b_2}$. 
Moreover, for a fixed pair $(a_1,b_1)$, $(a_2,b_2)$ of this kind,
the number of values $c$ for which $(a_1,b_1)$ and $(a_2,b_2)$ 
both belong to $\Pi_c$ is at most the constant degree $d$ of $F$,
unless $F$ vanishes identically on the two ``vertical'' lines
$(a_1,b_1)\times\reals$, $(a_2,b_2)\times\reals$, an assumption
that we adopt for our analysis.
%\micha{So we have to assume this property!}

It then follows, using the Cauchy-Schwarz inequality, that 
$$
I \ge \frac1d \sum_{c\in C} M_c^2 \ge \frac{\left(\sum_{c\in C} M_c\right)^2}{dn} 
= \frac{M^2}{dn} .
$$
The next step of the analysis is to derive an upper bound on $I$. 
On one hand this is an instance of a fairly standard point-curve incidence 
problem, which can be tackled using well established machinery, such as
the incidence bound of Pach and Sharir~\cite{PS98}, or, more fundamentally,
the crossing-lemma technique of Sz\'ekely~\cite{Sz97} (on which the analysis
in \cite{PS98} is based). However, to apply this machinery, there are several 
issues that need to be addressed: (a) The curves of $\Gamma$ are not necessarily 
distinct, or, more generally, many pairs of them might (partially) overlap,
at common irreducible components.
(b) We need to bound the number of intersections of any pair of distinct curves,
and to bound the number of curves that pass through any pair of points of $\Pi$.

As it turns out, and somewhat surprisingly, the first issue is the major hurdle
in the handling of the general problem. Before expanding upon this point,
let us first see how the technique continues in the ideal situation where 
all the curves of $\Gamma$ are distinct and non-overlapping. In this case
we have $|\Pi|=n^2$ distinct points and $|\Gamma|=n^2$ distinct curves in the $xy$-plane.
Using standard algebraic considerations, one can show that each pair of curves
intersect in $O(1)$ points (assuming the degree $d$ of $F$ to be constant),
and each pair of points can be incident to at most $O(1)$ common curves.
In this case, the techniques of \cite{PS98,Sz97} can be applied to yield
$$
I(\Pi,\Gamma) = O\left( |\Pi|^{2/3}|\Gamma|^{2/3} + |\Pi| + |\Gamma| \right)
= O(n^{8/3}) .
$$
Combining this with the lower bound on $I$, we get $M^2/n = O(n^{8/3})$,
or $M=O(n^{11/6})$.

Note that this improves considerably the bound $n^{1.95}$ in \cite{ER00,ESz}.

Let us return to the issue of coincidence or overlapping of the curves in $\Gamma$.
In the special instance that we study in this paper we use a concrete ad-hoc
argument that exploits the special geometric and algebraic structure of the
specific problem, allowing us to control the amount of coincidences and 
overlapping of curves.  What we are still missing, for the general problem,
is a general argument that if there are many coincident or overlapping pairs 
of curves, then $F$, $A$, $B$, and $C$ must have a special structure, similar 
to those established in \cite{ER00,ESz}. We find it rather strange that such 
a special structure (or the lack thereof) is manifested in the coincidence or 
overlapping (or the lack thereof) of the curves of $\Gamma$, and would like 
to better understand this connection.

An additional discussion of these and related issues is given in the 
concluding section.

%%%%%%%%%%%%%%%%%%%%%%%%%%%%%%%%%%%%%%%%%%%%%%%%%%%%%%%%%%%%%%%%%%%%%%%%%%%%%%%%

\section{Distinct distances from three points in the plane} \label{sec:3pts}

We recall the problem:
Let $p_1,p_2,p_3$ be three points in the plane, and let $P$ be a set of $n$ other 
points in the plane. The goal is to obtain a lower bound for the number of 
distinct distances between $p_1$, $p_2$, $p_3$ and the points of $P$, when
$p_1$, $p_2$, and $p_3$ are noncollinear.

\begin{figure}[htb]
\begin{center}
\input{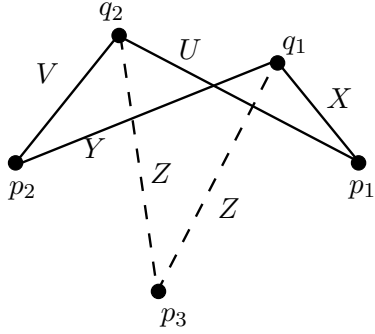}
\caption{A configuration involving the three fixed points $p_1,p_2,p_3$ and two
other points $q_1,q_2$ at equal distances from $p_3$. 
(The symbols $X,Y$, etc. are the squares of the lengths of the respective segments.)}
\label{XYUV}
\end{center}
\end{figure}

We may assume, without loss of generality, that $p_1=(1,0)$, $p_2=(-1,0)$, and $p_3=(a,b)$,
for $b\ne 0$. For a pair of points $q_1=(x,y)$ and $q_2=(u,v)$, we denote the 
squares of their distances from $p_1$ and $p_2$ as
\begin{align}\label{TwoDist}
X & = |p_1q_1|^2 = (x-1)^2+y^2 \nonumber \\
Y & = |p_2q_1|^2 = (x+1)^2+y^2 \\
U & = |p_1q_2|^2 = (u-1)^2+v^2 \nonumber \\
V & = |p_2q_2|^2 = (u+1)^2+v^2 \nonumber .
\end{align}
See Figure~\ref{XYUV}.

Let $P$ denote the set of the $n$ other given points.
We are going to estimate the number $Q$ of pairs $(q_1,q_2)\in P^2$, with $q_1\ne q_2$, 
which have equal distances from $p_3$. We will derive an upper bound 
and a lower bound for $Q$, and the comparison of these bounds will yield
the asserted lower bound on the number of distinct distances.

Before plunging into the analysis, we note that the problem at hand is
indeed a special instance of the general problem as reviewed in the introduction.
A point $q=(x,y)$ determines three squared distances $X=|p_1q|^2$, $Y=|p_2q|^2$,
and $Z=|p_3q|^2$ from the three respective points $p_1$, $p_2$, and $p_3$.
These distances must satisfy a polynomial equation $F(X,Y,Z)=0$; one can show
that this is a quadratic equation, although we will not make explicit use of this fact.
The $n$ points of $P$ determine $n$ triples $(X,Y,Z)$ at which $F$ vanishes.
If we denote by $\D$ the set of distinct distances between $p_1$, $p_2$, $p_3$ 
and the points of $P$, then $F$ vanishes at $n$ points of $\D\times\D\times\D$,
and our goal is in fact to obtain an upper bound for $n$ that is subquadratic
in $\kappa=|\D|$.

The analysis proceeds as follows.
For fixed values of $X$ and $V$, we define a planar curve $\gamma_{X,V}$,
in a parametric plane with coordinates $Y,U$,
which is the locus of all points $(Y,U)$ that, together with $X$ and $V$,
correspond to a pair of points 
$q_1=(x,y)$, $q_2=(u,v)$, so that these parameters satisfy (\ref{TwoDist}) and 
$|p_3q_1|=|p_3q_2|$, namely, 
$$
(x-a)^2+(y-b)^2 = (u-a)^2+(v-b)^2 ,
$$
or
\begin{equation}\label{ThirdDist}
x^2+y^2-2ax-2by = u^2+v^2-2au-2bv .
\end{equation} 
That is, $(Y,U)\in\gamma_{X,V}$ if and only if the following equations,
which result from a suitable combination of (\ref{TwoDist}) and (\ref{ThirdDist}),
have a common solution $(x,y,u,v)$.
\begin{align}\label{Hypert}
X & = (x-1)^2+y^2 \nonumber \\
V & = (u+1)^2+v^2 \nonumber \\
\tfrac12 (V-X) & = (1-a)x-by+(1+a)u+bv \\
Y & = X+4x \nonumber \\
U & = V-4u \nonumber.
\end{align}

Note that, given $X,Y,U,V$,
we can easily recover the corresponding coordinates $(x,y)$
and $(u,v)$, up to multiplicity of at most $4$, by observing 
that each of the two triangles $\triangle p_1p_2q_1$,
$\triangle p_1p_2q_2$, is fixed, up to a possible reflection 
about $p_1p_2$.  Algebraically, the coordinates
$x$ and $u$ are uniquely determined from the fourth and fifth
equations of (\ref{Hypert}), and the absolute values of $y$ and $v$
are then uniquely determined from the first two equations of
(\ref{Hypert}).

The third equation of (\ref{Hypert}) enforces the constraint that 
$(Y,U)\in \gamma_{X,V}$, and can be used to obtain the
algebraic equation of $\gamma_{X,V}$. That is, we have
\begin{align} \label{eqgamma}
\tfrac12 (V-X) & = \tfrac14(1-a)(Y-X) + \tfrac14(1+a)(V-U) \\
& + b\left(
\left(V-\tfrac14(V-U+4)^2\right)^{1/2} -
\left(X-\tfrac14(Y-X-4)^2\right)^{1/2} 
\right) , \nonumber
\end{align}
and we can turn this equation into a polynomial equation (in
$Y$ and $U$, regarding $X$ and $V$ as fixed parameters) of
degree four, as can be easily verified.

As remarked in the overview of the general problem in the 
introduction, a major technical hurdle that we need to overcome is
the possibility that many pairs of curves $\gamma_{X,V}$ coincide 
or overlap (in a common irreducible component).

For example, when $b=0$ (i.e., $p_1,p_2,p_3$ are collinear),
the equations (\ref{eqgamma}) are of parallel lines, all of the form
$U=\frac{1-a}{1+a}Y + c(X,V)$, where $c(X,V)$ is linear in $X$ and $V$.
In this case many curves can coincide with one another, and
the multiplicity of a curve can be as high as $\Theta(\kappa)$, 
where $\kappa$ is the number of distinct distances between 
$p_1$, $p_2$, $p_3$ and the points of $P$. As will be seen later,
this will cause our analysis to break down, in the sense that
in this case all we will be able to show is the trivial lower
bound $\kappa=\Omega(\sqrt{n})$. 
This will be further elaborated in a remark given at the 
end of the analysis.

Fortunately, as we next argue, when $b\ne 0$, the amount of coincidence 
or overlap between the curves is very limited.
More precisely, we have the following result.

\begin{proposition} \label{overlap} 
Each irreducible component of any curve $\gamma_{X,V}$ can be shared by at most
three other curves.
\end{proposition} 

\noindent{\bf Proof.}
We first observe that each curve $\gamma_{X,V}$ is bounded because, for $X,V$ fixed,
the point $q_1$ lies on a circle of radius $\sqrt{X}$ centered at $p_1$, and
the point $q_2$ lies on a circle of radius $\sqrt{V}$ centered at $p_2$.
This is easily seen to imply that any $(Y,U)\in \gamma_{X,V}$ must satisfy
$Y\le (2+\sqrt{X})^2$ and $U\le (2+\sqrt{V})^2$.

Let us consider an irreducible component $\gamma'_{X,V}$ of some curve $\gamma_{X,V}$
and a point on it, $(Y_0,U_0)$, such that $Y_0$ is maximal among all points 
of $\gamma'_{X,V}$. We will show that, given the point $(Y_0,U_0)$, the parameters
$X$ and $V$ can be recovered, up to multiplicity $4$.

Since $(U_0,Y_0)$ is $Y$-extremal, it has to satisfy the
equations $H(U_0,Y_0) = H_U(U_0,Y_0) = 0$, where $H=0$ is the algebraic 
equation of $\gamma_{X,V}$ given in (\ref{eqgamma}). The second equation is
$$
H_U(U_0,Y_0) = -\frac{1+a}4 + \frac{b}{4}\cdot
\frac{V-U_0+4}{(V-\tfrac14 (V-U_0+4)^2)^{1/2}} = 0 ,
$$
or
$$
V-\tfrac14 (V-U_0+4)^2) = 
\left( \frac{b}{1+a} \right)^2 (V-U_0+4)^2 .
$$
This is a quadratic equation in $V$ whose leading coefficient, namely
$\frac14 + \left(\frac{b}{1+a}\right)^2$,
is strictly positive, so it has at most two solutions.
(Note that the case $a=-1$, i.e., the case where $p_2$ and $p_3$
are co-vertical, is special, and yields the single solution
$V=U_0-4$.)

Next, fixing $V$ to be one of these two roots, the equation 
(\ref{eqgamma}) becomes an equation in $X$ of the form
$$
L(X) = (X-\tfrac14 (Y_0-X-4)^2)^{1/2} ,
$$
where $L(X)$ is a linear expression in $X$. Squaring this,
we obtain a quadratic equation in $X$ whose leading coefficient 
is strictly positive. Again, we obtain at most two solutions for 
$X$, for each value of $V$, for a total of at most four pairs $(X,V)$.

To sum up, we have shown that each irreducible 
curve $\gamma'$ can be a component of at most four curves $\gamma_{X,V}$,
as asserted.
$\Box$

%%%%%%%%%%%%%%%%%%%%%%%%%%%%%%%%%%%%%%%%%%%%%%%%%%%%%%%%%%%%%%%%%%%%%%%%%%

We continue with the analysis of the Erd{\H o}s--Lov\'asz--Vesztergombi problem.
Let $\D$ denote the set of distinct squared distances between $p_1,p_2,p_3$
and the points of $P$, and put $\kappa=|\D|$. 
Let $\Gamma$ denote the set of all curves $\gamma_{X,V}$, for $X,V\in \D$.

Every ordered pair $(q_1,q_2)$ of distinct points of $P$ with $|p_3q_1| = |p_3q_2|$
generate a quadruple $(X,Y,U,V)\in\D^4$, where
$$
X = |p_1q_1|^2,\qquad
Y = |p_2q_1|^2,\qquad
U = |p_1q_2|^2,\qquad
V = |p_2q_2|^2,
$$
such that $(Y,U)\in\gamma_{X,V}$.

The number $Q$ of these pairs $(q_1,q_2)$, introduced earlier, is proportional 
to the number of such quadruples (or incidences), because, as argued earlier,
each quadruple $(X,Y,U,V)$ can arise from at most four pairs $(q_1,q_2)$. 
We obtain a lower bound for $Q$, in complete analogy to the approach 
sketched in the introduction, as follows.  For each $Z\in\D$, denote by 
$P_Z$ the set of points at squared distance $Z$ from $p_3$.
Then, using the Cauchy-Schwarz inequality, we obtain
\begin{equation} \label{eqQ}
Q = \sum_{Z\in\D} \binom{|P_Z|}{2} = \frac12 \sum_{Z\in\D} |P_Z|^2
- \frac12 \sum_{Z\in\D} |P_Z| \ge
\frac{1}{2\kappa} \left( \sum_{Z\in\D} |P_Z| \right)^2 - \frac{n}{2} =
\frac{n^2}{2\kappa} - \frac{n}{2} .
\end{equation} 
To obtain an upper bound for $Q$, we bound the number of incidences between
the curves $\gamma_{X,V}$ and the points $(Y,U)$.  For this, we
apply Sz\'ekely's technique~\cite{Sz97}, which is based 
on the crossing lemma. This is also the approach used in
the proof of the incidence bound in Pach and Sharir~\cite{PS98},
but the possible overlap of curves, both in the primal 
and in the dual settings (see below for details), requires some 
extra (and more explicit) care in the application of the technique. 

In more detail, denote by $\Pi$ the set $\D^2$ of the $\kappa^2$ 
points $(Y,U)$, and let $\Gamma$ denote the (possibly multi-)set 
of the curves $\gamma_{X,V}$.
We begin by constructing a plane embedding 
of a multigraph $G$, whose vertices are the points of $\Pi$, 
and each of whose edges connects a pair $\pi_1=(Y_1,U_1)$,
$\pi_2=(Y_2,U_2)$ of points that lie on the same curve $\gamma_{X,V}$ 
and are consecutive along (some connected component of)
$\gamma_{X,V}$; one edge for each such curve (connecting $\pi_1$
and $\pi_2$) is generated.

%\micha{The next few paragraphs should be read VERY carefully.
%I tried several other arguments, and each one of them has crashed.
%Maybe there is a simpler argument, but if you think of one, 
%triple check it before you claim victory.}

A major potential problem with this construction is that
the edge multiplicity in $G$ may not be bounded (by a constant).
More concretely, we want to avoid edges $(\pi_1,\pi_2)$ whose 
multiplicity exceeds $16$.
We pass to a dual parametric plane, in which the roles of $(X,V)$
and $(Y,U)$ are interchanged, so points $(Y,U)$ of $\Pi$ become
dual curves that we denote as $\gamma_{Y,U}^*$, and curves 
$\gamma_{X,V}$ become dual points $(X,V)$.
By the symmetric nature of the definition, we have
$(Y,U)\in \gamma_{X,V}$ if and only if $(X,V)\in \gamma_{Y,U}^*$.
Hence, if the multiplicity of the edge connecting $(Y_1,U_1)$ and 
$(Y_2,U_2)$ is larger than $16$ then the dual curves 
$\gamma_{Y_1,U_1}^*$ and $\gamma_{Y_2,U_2}^*$ intersect in 
more than $16$ points, and therefore, since each is the zero 
set of a polynomial of degree $4$, B\'ezout's theorem implies 
that they must overlap in a common irreducible component.

Note that, given $(Y_1,U_1)$, the dual curve $\gamma_{Y_1,U_1}^*$,
having degree $4$, has 
at most four irreducible components, and, by Proposition~\ref{overlap},
applied in the dual plane, each such component can be shared by at 
most three other dual curves. That is, each $(Y_1,U_1)$ has at most 
$12$ ``problematic'' neighbors that we do not want to connect it to;
for any other point, the multiplicity of the edge connecting 
$(Y_1,U_1)$ with that point is at most $16$; more precisely, at
most $16$ curves $\gamma_{X,V}$ pass through both points.

Consider a point $(Y_1,U_1)$ and one of its bad neighbors
$(Y_2,U_2)$; that is, they are consecutive points along many curves.
Let $\gamma_{X,V}$ be one of the curves along which
$(Y_1,U_1)$ and $(Y_2,U_2)$ are neighbors. Then, rather than connecting
$(Y_1,U_1)$ to $(Y_2,U_2)$ along $\gamma_{X,V}$, we continue along 
the curve past $(Y_2,U_2)$ until we reach a good point for 
$(Y_1,U_1)$, and then connect $(Y_1,U_1)$ to that point (along 
$\gamma_{X,V}$).  We skip over at most $12$ points in the process, 
but now, having applied this ``stretching'' to each pair of bad 
neighbors, each of the modified edges has multiplicity at most $16$.

The number of new edges in $G$ is at least 
$I(\Pi,\Gamma)-c|\Gamma|$, for a suitable constant $c$,
where the term $c|\Gamma|$ accounts for the number of connected 
components of the curves---for components with fewer than $14$
incident points, there might be no edge drawn along that component. 

The final ingredient needed for this technique is an upper bound
on the number of crossings between the edges of $G$.  Each such 
crossing is a crossing between two curves of $\Gamma$.  Even 
though the two curves might overlap in a common irreducible 
component (where they have infinitely many intersection points, 
none of which is a crossing), the number of proper crossings between 
them is $O(1)$, as follows, for example, from the Milnor--Thom and 
B\'ezout's theorems. Finally, because of the way the drawn edges have 
been stretched, the edges now may overlap one another, and then
a crossing between two curves may be claimed by more than one pair
of edges.  Nevertheless, since no edge straddles through more than 
$12$ points, the number of pairs that claims a specific crossing
is $O(1)$.  Hence, we conclude that the total number of edge 
crossings in $G$ is $O(|\Gamma|^2)$.

We can now continue by applying the crossing lemma, exactly
as done in many earlier works (e.g., see \cite{PS98,Sz97}), and 
conclude that
$$
I(\Pi,\gamma) = O\left( |\Pi|^{2/3}|\Gamma|^{2/3} + 
|\Pi| + |\Gamma| \right) .
$$
Since $|\Pi| = |\Gamma| = O(\kappa^2)$, it follows that
%%%%%%%%%%%%%%%%%%%%%%%%%%%%%%%%%%%%%%%%%%%%%%%%%%%%%%%%%
$$
Q = O(I(\Pi,\Gamma')) = O(\kappa^{8/3}) .
$$
Comparing this with the lower bound in (\ref{eqQ}), we obtain
$$
\frac{n^2}{2\kappa} - \frac{n}{2} = O(\kappa^{8/3}),\qquad\text{or}\qquad 
\kappa = \Omega(n^{6/11}) .
$$
That is, we have obtained the following main result of this section.
%%%%%%%%%%%%%%%%%%%%
\begin{theorem} \label{thm:elv}
The number of distinct distances between three noncollinear points and 
$n$ other points in the plane is $\Omega(n^{6/11})$.
\end{theorem}
%%%%%%%%%%%%%%%%%%%%

\noindent{\bf Remark.}
Returning to a comment made earlier, we note that the preceding analysis
can also be applied when $p_1,p_2,p_3$ are collinear. In this case the maximum
multiplicity of a curve $\gamma_{X,V}$ (which is a line of a fixed slope in this case)
is $\kappa$, because each value of $X$ determines a unique value of $V$ that yields
the same curve (that is, line). 
We can carry out the analysis by considering the worst case, where we have
only $O(\kappa)$ distinct curves, each with multiplicity $\kappa$. 
In this case the upper bound on $Q$ is
$$O(\kappa\cdot (\kappa^{2/3}(\kappa^2)^{2/3} + \kappa + \kappa^2)) = O(\kappa^3) ,$$
so $n^2/\kappa = O(\kappa^3)$, or $\kappa=\Omega(\sqrt{n})$. This lower bound matches
the upper bound in the construction of Elekes and Szab\'o~\cite{ESz}, but it is
totally trivial, because if there were fewer than $\frac12\sqrt{n}$ distinct 
distances, $P$ would have contained fewer than $n$ points. The present remark 
is made in order to highlight the significance of the non-overlapping of the curves.

%%%%%%%%%%%%%%%%%%%%%%%%%%%%%%%%%%%%%%%%%%%%%%%%%%%%%%%%%%%%%%%%%%%%%%%%%%%%%%

\section{Discussion}

We have studied an interesting problem in combinatorial geometry, 
concerning the number of distinct distances between three noncollinear 
points and $n$ other points in the plane, which can also be regarded as 
a special instance of a more general problem, concerning the number of points
of a triple Cartesian product at which a given trivariate polynomial can
vanish. The general problem has been tackled in \cite{ER00,ESz}, but we
have bypassed this general approach, replaced it with a different novel
general approach, and combined it with a direct ad-hoc 
technique, and have thereby managed to improve (a) the bound
on the number of zeros of the specific polynomial $F$ that arises
in our setup, and (b) the earlier bound for the specific distinct 
distances problem, as obtained in \cite{ESSz}.

Certainly, one of the main open problems is to understand better the
structure of the general problem. In particular, what is the connection
between low multiplicities of the curves that we define and the structure
of the polynomial $F$ (as provided in \cite{ER00,ESz})? A more concrete
formulation of the problem is to find a general technique for showing that
if our curves have high multiplicity then $F$ must have the special form
given in \cite{ER00,ESz} (or perhaps some other special form?)

In parallel, it would be interesting to identify other special instances
of the general problem, and apply our machinery to obtain new or improved 
bounds for them. One such instance, that we are currently studying, is the 
following problem, studied by Elekes, Simonovits and Szab\'o~\cite{ESSz}. 
Let $p_1,p_2,p_3$ be three distinct points in the plane, and, for $i=1,2,3$, let 
$\C_i$ be a family of $n$ unit circles that pass through $p_i$. The goal is to
obtain a subquadratic upper bound on the number of \emph{triple points}, 
which are points that are incident to a circle of each family. 
Elekes et al.~\cite{ESSz} have shown that the number of such points 
is $O(n^{2-\eta})$, for some constant parameter $\eta>0$ (that they did 
not specify), as an application of a more general technique that they have 
developed (see also other references in \cite{ESSz}). 

This problem too fits into the general framework, and if the multiplicities 
of the resulting curves could be shown to be under control, we would obtain
that the number of triple points is $O(n^{11/6})$, improving the bound and 
making it more concrete.

Another such possible problem is the following. Let $P_1$, $P_2$, and $P_3$
be three sets of $n$ points each, so that each set $P_i$ is contained 
in some line $\ell_i$, for $i=1,2,3$. How many unit-area triangles are
determined by triples of points in $P_1\times P_2\times P_3$?

It turns out that our technique can also be applied to the following problem. 
Let $\gamma$ be a small-degree algebraic curve in the plane, and let $P$ be a
set of $n$ points lying on $\gamma$. How many distinct distances must there 
be between the points of $P$? We have recently obtained \cite{SSS} an 
improved lower bound of $\Omega(n^{4/3})$ for a special bipartite version 
of this problem, where we have two sets $P_1$, $P_2$ of $n$ points each, 
lying on two respective lines in the plane, which are neither parallel nor 
orthogonal. Later, and very recently, Pach and de Zeeuw~\cite{PdZ} have
extended the machinery to the general case, with a similar lower bound.

We note that the application of our technique to this problem is interesting 
because it does not seem to fit into the paradigm of a polynomial vanishing 
on a 3-dimensional Cartesian product, but it nevertheless benefits from our approach.

In conclusion, we note that the bounds that we have obtained are asymmetric.
In the notation of the general problem, if the sets $A$, $B$, $C$ are of 
different sizes, our bound on $|Z(F)|$ becomes $O(|A|^{2/3}|B|^{2/3}|C|^{1/2})$,
with similar asymmetric consequences for the two specific problems.
However, the general problem is fully symmetric in $A$, $B$, and $C$, so one
would definitely expect a bound that is symmetric in the sizes of the 
three sets. We leave this refinement of the bound as an open problem.


\begin{thebibliography}{}

\bibitem{El}
G. Elekes, 
Sums versus products in number theory, algebra and Erd{\H o}s geometry--A survey,
in {\it Paul Erd{\H o}s and his Mathematics} II,
Bolyai Math. Soc., Stud. 11, Budapest, 2002, pp.~241--290.

\bibitem{ESSz}
G. Elekes, M. Simonovits and E. Szab\'o,
A combinatorial distinction between unit circles and straight lines:
How many coincidences can they have?
{\it Combinat. Probab. Comput.} 18 (2009), 691--705.

\bibitem{ER00}
G.\ Elekes and L.\ R\'onyai, 
A combinatorial problem on polynomials and rational functions,
{\it J.\ Combinat.\ Theory Ser.\ A}, 89 (2000), 1--20.

\bibitem{ESz}
G. Elekes and E. Szab\'o,
How to find groups? (and how to use them in Erd{\H o}s geometry?),
{\it Combinatorica} 32 (2012), 537--571.

\bibitem{ELV}
P. Erd{\H o}s, L. Lov\'asz, and K. Vesztergombi,
On the graph of large distance,
{\it Discrete Comput. Geom.} 4 (1989), 541--549.

\bibitem{PS98}
J. Pach and M. Sharir,
On the number of incidences between points and curves,
{\it Combinat. Probab. Comput.} 7 (1998), 121--127.

\bibitem{PdZ}
J. Pach and F. de Zeeuw,
Distinct distances on algebraic curves in the plane,
in arXiv:1308.0177.

\bibitem{SSS}
M. Sharir, A. Sheffer, and J. Solymosi,
Distinct distances on two lines,
{\it J. Combinat. Theory}, Ser. A 20 (2013), 1732--1736.

\bibitem{Sz97}
L. Sz\'ekely,
Crossing numbers and hard Erd{\H o}s problems in discrete geometry,
{\it Combinat. Probab. Comput.} 6 (1997), 353--358.

\end{thebibliography}
\end{document}